\documentclass[11pt,reqno]{amsart}

\usepackage{amsmath}

\usepackage{amsthm}

\usepackage{amssymb}

\usepackage{graphics}

\theoremstyle{plain}

\newtheorem{theorem}{Theorem}[section]

\newtheorem{proposition}[theorem]{Proposition}

\newtheorem{corollary}[theorem]{Corollary}

\theoremstyle{definition}

\newtheorem{examples}[theorem]{Examples}

\theoremstyle{remark}

\newtheorem{remark}[theorem]{Remark}
\newtheorem{remarks}[theorem]{Remarks}

\DeclareMathOperator{\zero}{\mathcal{Z}}

\newcommand{\Ew}{\mathbb{E}} 
\newcommand{\PP}{\mathbb{P}} 
\newcommand{\holder}{H\"{o}lder }

\begin{document}

\title[Cantor set zeros of Brownian motion -- Cantor function]{Cantor set zeros of one-dimensional Brownian motion minus Cantor function}

\author[J. Ruscher]{Julia Ruscher}

\address{
Fachbereich Mathematik\\
Technische Universit\"{a}t Berlin\\
Strasse des 17. Juni 136\\
D-10623 Berlin\\
and
\newline
Institute for Advanced Study\\
Einstein Drive\\
Princeton, NJ 08540
}

\email{ruscher@math.tu-berlin.de}




\begin{abstract}
In \cite{ABPR} it was shown by Antunovi\'{c}, Burdzy, Peres, and Ruscher that a Cantor function added to one-dimensional Brownian motion has zeros in the middle $\alpha$-Cantor set, $\alpha \in (0,1)$, with positive probability if and only if $\alpha \neq 1/2$.
We give a refined picture by considering a generalized version of middle $1/2$-Cantor sets. By allowing the middle $1/2$ intervals to vary in size around the value $1/2$ at each iteration step we will see that there is a big class of generalized Cantor functions such that if these are added to one-dimensional Brownian motion, there are no zeros lying in the corresponding Cantor set almost surely.

\end{abstract}
\maketitle

\section{Introduction}

Let $B(t)$ be standard one-dimensional Brownian motion with $B(0)=0$, let $C_\alpha$ be be a middle $\alpha$-Cantor set for $\alpha \in (0,1)$ on the interval $[1,2]$ and let $f_\alpha \colon [1,2] \to [0,1]$ be the corresponding Cantor function.
For any function $g$ defined on some interval $I \subset \mathbb{R}^+$ denote by $\zero(g)$ the set of zeros of $g$ in $(0,\infty)$.

Taylor and Watson (Example 3 in \cite{TaylorWatson}) showed that Brownian does not intersect the graph of the middle $1/2$-Cantor function restricted to the Cantor set, even though the projection of this set on the vertical axis is an interval.
Surprisingly, the middle $1/2$-Cantor function is an exception in the following sense.

\begin{theorem}[\cite{ABPR}] \label{thm: ABPR}
$\mathbb{P}(\zero(B-f_\alpha) \cap C_\alpha \neq \emptyset)>0$ holds if and only if $\alpha \neq 1/2$.
\end{theorem}
In the paper \cite{ABPR} part of this result is used to prove that for every $\beta<1/2$ there is a $\beta$-\holder continuous function $f \colon \mathbb{R}^+ \to \mathbb{R}$ such that the set $\zero(B-f)$ has isolated points with positive probability.

Theorem \ref{thm: ABPR} is the main motivation of our present work. We investigate a more general class of Cantor functions defined in section \ref{section: Cantor function def}. The $n$-th approximation of the middle $1/2$-Cantor function increases on some intervals of length $4^{-n}$ and is constant elsewhere. We now allow these intervals to vary in length for every iteration step meaning that for a positive, real sequence $(a_n)$ the length of intervals where the function increases is $4^{-n}a_n^2$ at iteration level $n$. See Figure \ref{fig:cantorII} for example.
Theorems \ref{thm:cantorset-zeros2} and \ref{thm: cantorset-zeros 2nd part} give conditions for which of these generalized Cantor functions $f_b$ the process $B-f_b$ has zeros in the corresponding generalized Cantor set with positive probability or with probability $0$, respectively.

Other questions concerning path properties of Brownian motion with variable drift lying outside the Cameron-Martin space have been studied recently, e.g. in \cite{ABPR}, \cite{APV}, \cite{PS} and \cite{R12}.

\section{Generalized Cantor sets and generalized Cantor function}\label{section: Cantor function def}

We now define a generalized Cantor set analogously to the standard Cantor set (see for instance section 3 of \cite{ABPR}). 

For a given positive, real sequence $(a_n)$ we define a sequence $(b_{n})$ by $b_{n}:= (2^{-n} \cdot a_{n})^2$. For this sequence we define a corresponding Cantor-type set and denote it by $C_{b}$.
Take a closed interval $I$ of length $b_0=a_0$. Let $\mathfrak{C}_{b_{1}}$ be the set consisting of two disjoint closed subintervals of $I$ of length $b_{1}$, the left one (for which the left endpoint coincides with the left endpoint of $I$) and the right one (for which the right endpoint coincides with the right endpoint of $I$).
Now continue recursively, if $J \in \mathfrak{C}_{b_{n}}$, then include in the set $\mathfrak{C}_{b_{n +1}}$ its left and right closed subintervals of length $b_{n +1}$. We define the set $C_{b_{n}}$ as the union of all the intervals from $\mathfrak{C}_{b_{n}}$.
For any $n$, the family $\mathfrak{C}_{b_{n}}$ is the set of all connected components of the set $C_{b_{n}}$.
The generalized Cantor set is a compact set defined as $C_b=\bigcap_{n=1}^\infty C_{b_{n}}$.

Now we construct a Cantor-type function corresponding to the generalized Cantor set above. Define the function $f_{b_{1}}$ so that it has values $0$ and $1$ at the left and the right endpoint of the interval $I$, respectively, value $1/2$ on $I\backslash {C}_{b_{1}}$ and interpolate linearly on the intervals in $\mathfrak{C}_{b_{1}}$. Recursively, construct the function $f_{b_{n +1}}$ so that for every interval $J =[s,t]\in \mathfrak{C}_{b_{n}}$, the function $f_{b_{n +1}}$ agrees with $f_{b_n}$ at $s$ and $t$, it has value $(f_{b_n}(s)+f_{b_n}(t))/2$ on $J \backslash C_{b_{n +1}}$ and interpolate linearly on the intervals in $\mathfrak{C}_{b_{n +1}}$. 

\begin{figure}

\includegraphics{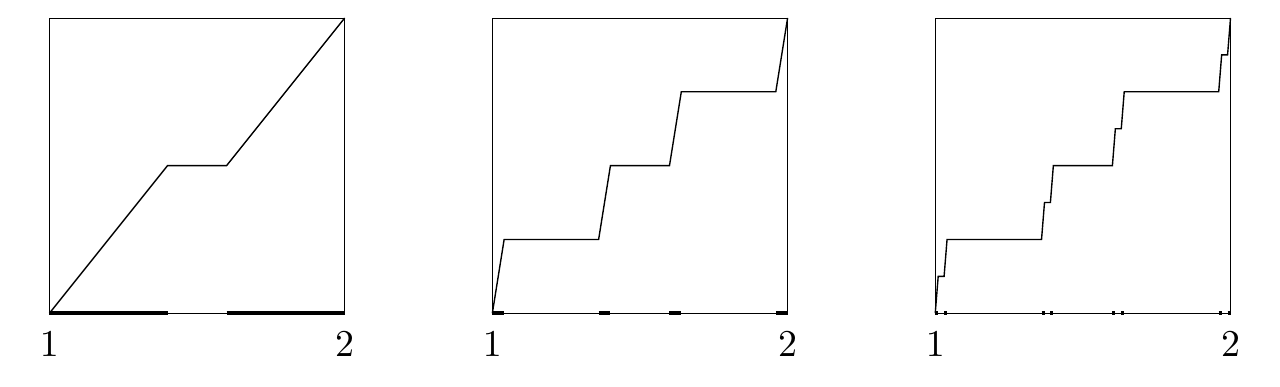}

\caption{First three approximations of the generalized Cantor function on the interval $[1,2]$ (functions $f_{b_n}$ from the construction for $b_1=0.4$, $b_2=0.04$ and $b_3=0.01$, or $a_1=\sqrt{1.6}$, $a_2=\sqrt{0.64}$ and $a_3=0.8$, respectively. Approximations of the generalized Cantor set $C_{b}$ are drawn in bold).}

\label{fig:cantorII}

\end{figure}

The sequence of functions $(f_{b_n})$ converges uniformly on $I$.
We define the generalized Cantor function $f_b$ as the limit $f_b=\lim_n f_{b_n}$.
For any $n$ and all $m \leq n$ the functions $f_b$ and $f_{b_n}$ agree at the endpoints of intervals $J \in \mathfrak{C}_{b_m}$. See Figure \ref{fig:cantorII} for an example.

Further, we fix an arbitrary small $\epsilon > 0$ and require the sequence $(a_n)$ to fulfill the condition
$a_n^2 - \frac{1}{2}a_{n+1}^2 \geq \epsilon a_n^2$ (or equivalently $b_n- 2b_{n+1} \geq \epsilon b_n$) for all $n$. Corollary \ref{rm: condition on a_n's} addresses the more general case of having a weaker condition $b_n- 2b_{n+1} > 0$ for all $n$. So we will only consider non-degenerate Cantor sets/ functions. Note that if the sequence $(a_n)$ is monotone decreasing, then $a_n^2 - \frac{1}{2}a_{n+1}^2 \geq \frac{1}{2} a_n^2$ holds for all $n$.

For simplicity, we will assume the initial interval $I$ to be $[1,2]$, and we can extend the function $f_{b}$ to $\mathbb{R}^+$, for instance by value $0$ on $[0,1)$ and value $1$ on $(2,\infty)$. 

\section{Zeros in the generalized Cantor set}\label{Zeros in the Cantor set}

The following theorem gives a condition for which sequences $(b_n)$ (or $(a_n)$ respectively) $B-f_b$ has zeros in the generalized Cantor set with positive probability.

\begin{theorem}\label{thm:cantorset-zeros2}
$\mathbb{P}(\zero(B-f_b) \cap C_b \neq \emptyset)>0$ holds if either
\\$\sum_{n=1}^{\infty} a_n < \infty$
, or
\\$\sum_{n=1}^{\infty} \frac{1}{a_n} < \infty$
.
\end{theorem}

\begin{remarks}
(i) Note that for geometric series $(a_n)$ the result was already shown in Theorem \ref{thm: ABPR} (i.e. Theorem 1.3. of \cite{ABPR}). For $a_n = \frac{1}{x^{n}}$ with some $x>1$ gives $b_{n}:= (2^{-n} \cdot \frac{1}{x^{n}})^2 = (2x)^{-2n}$ which corresponds to $\alpha = 1- \frac{1}{x}$ in Theorem \ref{thm: ABPR}. So for $\sum_{n=1}^{\infty} a_n < \infty$ Theorem \ref{thm:cantorset-zeros2} extends Theorem \ref{thm: ABPR} for convergent series $(\sum a_n)$ that increase slower than geometric series, for instance take $a_n = \frac{1}{n^d}$ with $d>1$.

(ii)The generalized Cantor function $f_b$ is not necessarily $\alpha$-\holder continuous for some $\alpha \in (0,1)$. If the sequence $a_n$ fulfills that $-\log 2 \cdot \frac{n}{\log a_n}$ converges to some $\sigma \in (0,1)$, then the corresponding generalized Cantor function is $\sigma$-\holder continuous. For example, the generalized Cantor function corresponding to the sequence $(n^{-d})$ with $d>1$ does not satisfy a \holder condition of any positive order.

\end{remarks}

To prove Theorem \ref{thm:cantorset-zeros2} we will use essentially the same method as in \cite{ABPR}.

\begin{proof}[Proof of Theorem \ref{thm:cantorset-zeros2}]
For an interval $I=[r,s] \in \mathfrak{C}_{b_n}$, define $Z_n(I)$ as the event $B(s) \in [f_b(r),f_b(s)]$, and the random variable $Z_{b_n} = \sum_{I \in \mathfrak{C}_{b_n}}\mathbf{1}(Z_n(I))$, where $\mathbf{1}(Z_n(I))$ is the indicator function of the event $Z_n(I)$.

Note that, by simple bounds on the transition density of Brownian motion, there is a constant $c_1>0$, such that for any sequence $b_{n}:= (2^{-n} \cdot a_{n})^2$, 

\begin{equation}\label{eq: hitting probabilities}
c_12^{-n} \leq \mathbb{P}(Z_n(I)) \leq 2^{-n} \ \text{ and } \ c_1 \leq \mathbb{E}(Z_{b_n}) \leq 1.
\end{equation}

If the event $Z_{b_n}>0$ happens for infinitely many $n$'s, then we can find a sequence of intervals $I_k=[r_k,s_k] \in \mathfrak{C}_{b_{n_k}}$, such that $f_{b}(r_k) \leq B(s_k) \leq f_{b}(s_k)$, thus $|B(s_k) - f_b(s_k)| \leq 2^{-n_k}$. Since $s_k \in C_{b_{n_k}}$, the sequence $(s_k)$ will have a subsequence converging to some $s \in C_b$, which satisfies $B(s)=f_b(s)$.
Therefore
\begin{equation}\label{eq: approximation of probabilities_1}
\mathbb{P}(\zero(B-f_b) \cap C_{b} \neq \emptyset) \geq \mathbb{P}(\limsup_{n\rightarrow \infty} \{Z_{b_n}>0\}).
\end{equation}
To bound the probabilities $\mathbb{P}(Z_{b_n}>0)$ from below we apply the Paley-Zygmund inequality 
(see Lemma 3.23 in \cite{MP}):
\begin{equation*}\mathbb{P}(Z_{b_n}>0) \geq (\mathbb{E} Z_{b_n})^2/\mathbb{E} (Z_{b_n}^2).
\end{equation*}
Therefore, we have to bound the second moment $\mathbb{E}(Z_{b_n}^2)$ from above. We will use the following expression for the second moment
\begin{equation}\label{eq: second moment}
\mathbb{E}(Z_{b_n}^2) = 2\sum_{I,J \in \mathfrak{C}_{b_n} \atop I < J} \mathbb{P}(Z_n(I))\mathbb{P}(Z_n(J) \mid Z_n(I)) + \mathbb{E}(Z_{b_n}),
\end{equation}
where by $I<J$ we mean that the interval $I$ is located to the left of the interval $J$.
Now we fix $n$ and intervals $I=[s_1,t_1]$ and $J=[s_2,t_2]$ from $\mathfrak{C}_{b_n}$, so that $I<J$. Let $x_i=f_b(s_i)$ and $y_i=f_b(t_i)$, for $i=1,2$. By the Markov property and the scaling property of Brownian motion, the process 

\[\widetilde{B}(t)=(t_2-t_1)^{-1/2}\Big(B(t_1+(t_2-t_1) t)-B(t_1)\Big),\]
is again a Brownian motion, independent of $\mathcal{F}_{t_1}$, and thus independent of the event $Z_n(I) \in \mathcal{F}_{t_1}$.

The event $Z_n(J)$ happens when $\widetilde{B}(1) \in \overline{J}$ for the interval $\overline{J}=[(t_2-t_1)^{-1/2}(x_2-B(t_1)) , (t_2-t_1)^{-1/2}(y_2-B(t_1))]$ of length $(t_2-t_1)^{-1/2}2^{-n}$.

\begin{figure}

\includegraphics{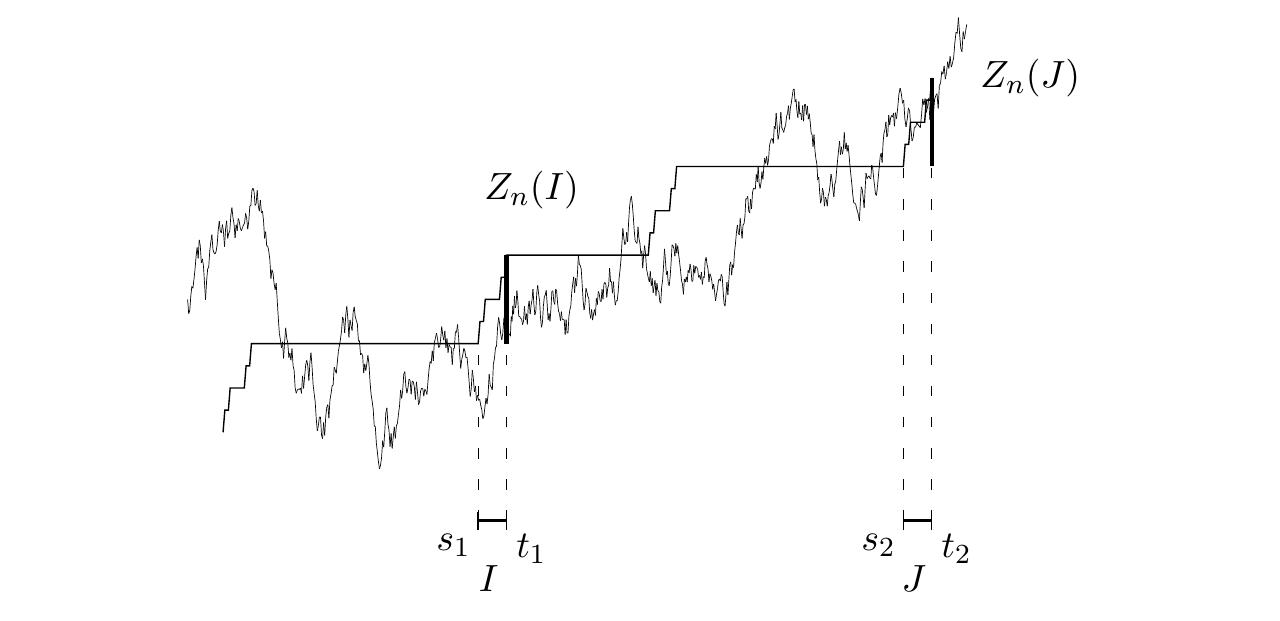}

\caption{Events $Z_n(I)$ and $Z_n(J)$ (graph of Brownian motion intersects two bold vertical intervals).}

\label{fig:intersection_with_cantorII}

\end{figure}

Fix intervals $I^0$ and $J^0$ in $\mathfrak{C}_{b_{\ell+1}}$, which are contained in a single interval in $\mathfrak{C}_{b_{\ell}}$. Assume $I^0 < J^0$ and label the intervals from $\mathfrak{C}_{b_n}$ contained in $I^0$ by $I_1, \dots, I_{2^{n-\ell-1}}$, and those contained in $J^0$ by $J_1, \dots, J_{2^{n-\ell-1}}$, so that $I_{i+1} < I_i$ and $J_j < J_{j+1}$. Set $I=I_i$ and $J=J_j$ for some $1 \leq i,j \leq 2^{n-\ell-1}$, and define $x_i$, $y_i$, $\widetilde{B}$ and $\overline{J}$ as before.
For $t_2-t_1$ we use the estimates $\epsilon b_\ell = \epsilon  2^{-2\ell}  a_\ell^2 \leq 2^{-2\ell}(a_\ell^2 - \frac{1}{2}a_{\ell+1}^2) = b_\ell - 2b_{\ell+1} \leq t_2-t_1 \leq b_\ell$. Conditional on $Z_n(I)$, the left endpoint of the interval $\overline{J}$ is at least $(x_2-y_1)(t_2-t_1)^{-1/2}$ and the right endpoint is at most $(y_2-x_1)(t_2-t_1)^{-1/2}$. Since $x_2-y_1 = (i+j-2)2^{-n}$ and $y_2-x_1 = (i+j)2^{-n}$ we have $\overline{J} \subset [(i+j-2)2^{-n}b_\ell^{-1/2}, (i+j)2^{-n}(\epsilon b_\ell)^{-1/2}]$. Because $\overline{J}$ has length at least $(b_\ell)^{-1/2}2^{-n} = \frac{2^{\ell-n}}{a_\ell}$ and at most $(\epsilon b_\ell)^{-1/2}2^{-n}= \frac{2^{\ell-n}}{\sqrt{\epsilon}a_\ell}$ we obtain

\begin{align*}
\mathbb{P}(Z_n(J_j) \mid Z_n(I_i)) &  = \mathbb{P}(\widetilde{B}(1) \in \overline{J} \mid Z_n(I))\\
& \leq \frac{2^{\ell-n}}{\sqrt{2\pi\epsilon}a_\ell}\exp\Big(-\frac{1}{2} {((j+i-2)\frac{2^{\ell-n}}{a_\ell})^2}\Big).
\end{align*}
By summing over $1 \leq i,j \leq 2^{n-\ell-1}$ it follows that

\begin{multline}\label{eq: conditional probabilities estimates}
\sum_{1 \leq i,j \leq 2^{n-\ell-1}}\mathbb{P}(Z_n(J_j) \mid Z_n(I_i))
\\
\leq 1 + \frac{1}{\sqrt{2\pi\epsilon}}
\sum_{k=1}^{\infty}(k+1)\frac{2^{\ell-n}}{a_\ell}\exp\Big(-\frac{1}{2}(\frac{2^{\ell-n}k}{a_\ell})^2\Big)
\end{multline}
where we used the trivial bound for $i=j=1$. The sum on the right hand side can be written as

\begin{equation*}
S(z)=z\sum_{k = 1}^\infty(k+1)\exp(-(kz)^2/2),
\end{equation*}
for $z=\frac{2^{\ell-n}}{a_\ell}$. Since $\exp(-t^2/2) \leq \exp(-t+1/2)$, we get that
\[
S(z) \leq e^{1/2}z\sum_{k \geq 1}(k+1)e^{-kz} = e^{1/2}z^{-1}\Big(\frac{z}{1-e^{-z}}\Big)^2(2e^{-z}-e^{-2z}).
\]
Since $z\mapsto (\frac{z}{1-e^{-z}})^2(2e^{-z}-e^{-2z})$ is a bounded function on $\mathbb{R}^+$, it follows by (\ref{eq: conditional probabilities estimates}) that for any fixed intervals $I^0$ and $J^0$ as above

\begin{equation}
\label{eq:second_moment_estimates_below_1/4_3}
\sum_{I,J \in \mathfrak{C}_{b_n} \atop I \subset I^0, J \subset J^0} \mathbb{P}(Z_n(J) \mid Z_n(I)) \leq 1+c_2 \frac{a_\ell}{2^{\ell-n}},
\end{equation}
for some constant $c_2>0$.

Therefore, summing the inequality in (\ref{eq:second_moment_estimates_below_1/4_3}) over all intervals $I^0$ and $J^0$ and $\ell=0, \dots , n-1$, and using it together with (\ref{eq: second moment}) and (\ref{eq: hitting probabilities}), we have
\[
\mathbb{E}(Z_{b_n}^2)
\leq 2^{-n+1} \sum_{\ell=0}^{n-1}\Big(2^\ell + 2^\ell c_2 \frac{a_\ell}{2^{\ell-n}}\Big) + 1 \leq 2\sum_{\ell=0}^n\Big(2^{-(n-\ell)} + c_2 a_\ell\Big).
\]
Now we see that $\mathbb{E}(Z_{b_n}^2)$ is bounded from above if $\sum_{n=1}^{\infty} a_n < \infty$. The lower bound in the second inequality in (\ref{eq: hitting probabilities}) and the Paley-Zygmund inequality imply that $\mathbb{P}(Z_{b_n}>0) \geq \mathbb{E}(Z_{b_n})^2/\mathbb{E}(Z_{b_n}^2)$ is bounded from below and the claim follows for the first case of the theorem from (\ref{eq: approximation of probabilities_1}).

Now for the second case pick intervals $I,J \in \mathfrak{C}_{b_n}$ such that $[s_1,t_1]=I < J =[s_2,t_2]$ and define $x_i$, $y_i$, $\widetilde{B}$ and $\overline{J}$ as before.
By $\ell$ denote the largest integer such that both $I$ and $J$ are contained in a single interval from $\mathfrak{C}_{b_\ell}$.
Assume that the event $Z_n(I)$ happens.
We see that the endpoints of the interval $\overline{J}$ are satisfying that

\[\frac{x_2-B(t_1)}{(t_2-t_1)^{1/2}} \geq 0 \ \text{ and } \  \frac{y_2-B(t_1)}{(t_2-t_1)^{1/2}} \leq 2^{-\ell}\cdot \sqrt{\frac{1}{\epsilon}}b_\ell^{-1/2} = \frac{1}{\sqrt{\epsilon} a_\ell}.\]

If $\frac{1}{a_\ell}$ is bounded, then the interval $\overline{J}$ is contained in a compact interval, which does not depend on the choice of $n$, $\ell$, $I$ or $J$. Using this and the fact that the length of $\overline{J}$ is bounded
with $ (b_\ell)^{-1/2}2^{-n} \leq |\overline{J}| \leq (\epsilon b_\ell)^{-1/2}2^{-n}$,  we get that for some positive constants $c_3$ and $c_4$ we have

\begin{equation}\label{eq: easy conditional probabilities}
c_3(b_\ell)^{-1/2}2^{-n} \leq \mathbb{P}(Z_n(J) \mid Z_n(I)) = \mathbb{P}(\widetilde{B}(1) \in \overline{J} \mid Z_n(I)) \leq c_4 (\epsilon b_\ell)^{-1/2}2^{-n}.
\end{equation}
Substituting  (\ref{eq: easy conditional probabilities}) and the upper bounds from (\ref{eq: hitting probabilities}) into (\ref{eq: second moment}), and summing over all intervals $I$ and $J$, we obtain

\[
\mathbb{E}(Z_{b_n}^2) \leq 1 + 2^{-n+1}\sum_{\ell=0}^{n-1}2^\ell2^{2(n-\ell-1)}c_4 (\epsilon b_\ell)^{-1/2}2^{-n} = 1 + c_5\sum_{\ell=0}^{n-1} \frac{1}{a_\ell}
\]
for some positive constant $c_5$.
We see that if $\sum_{n=1}^{\infty} \frac{1}{a_n} < \infty$, then $\mathbb{E}(Z_{b_n}^2)$ is bounded from above by a constant not depending on $n$, and the claim follows.

\end{proof}

The following theorem gives a condition for which sequences $(b_n)$ (or $(a_n)$ respectively) $B-f_b$ has no zeros in the generalized Cantor set almost surely.


\begin{theorem}\label{thm: cantorset-zeros 2nd part}
(i) If there is a sequence $(c_n)$ with $\sum_{n=1}^{\infty} c_n = \infty$ and a fixed but arbitrary small $\delta > 0$  and an $n_0>0$ such that for all $n \geq n_0$
\begin{align}\label{inequality of thm: cantorset-zeros 2nd part}
\frac{1}{c_n}\geq a_n \geq \sqrt{\frac{1}{(8-\delta)\epsilon\ln{\frac{1}{c_n}}}},
\end{align}
then $\mathbb{P}(\zero(B-f_b) \cap C_b \neq \emptyset)=0$.

(ii) If $\sum_{\ell = 1}^n \frac{1}{a_\ell} \rightarrow \infty$ for $n \rightarrow \infty$ and if there is a constant $C>0$ such that $\sum_{\ell = 1}^n \frac{1}{a_\ell} \leq C \log n$ for all $n>0$, then $\mathbb{P}(\zero(B-f_b) \cap C_b \neq \emptyset)=0$.

\end{theorem}

\begin{examples}
The sequence $(a_n)$ defined by $a_n = \sqrt{\frac{1}{\ln{n}}}$ fulfills the conditions of the Theorem \ref{thm: cantorset-zeros 2nd part}(i) (to see that choose $c_n= \frac{1}{n}$) and the sequence $(a_n)$ defined by $a_n = n$ fulfills the conditions of the Theorem \ref{thm: cantorset-zeros 2nd part}(ii).

By Theorem \ref{thm: cantorset-zeros 2nd part}(i) also applies to sequences $(a_n)$ where every element of the sequence is chosen from a fixed finite set of a numbers. Therefore, we see that $\mathbb{P}(\zero(B-f_b) \cap C_b \neq \emptyset)=0$ holds for all these sequences.
\end{examples}

Note that, if two sequences $(a_n)$ and $(a'_n)$ only differ by finitely many numbers, and if one of the sequences fulfills the conditions of one of the Theorems \ref{thm:cantorset-zeros2} or \ref{thm: cantorset-zeros 2nd part}, then the other sequence fulfills the conditions of the same theorem.

\begin{proof}[Proof of Theorem \ref{thm: cantorset-zeros 2nd part}]
For an interval $I\in \mathfrak{C}_{b_n}$, define $Y_n(I)$ as the event that Brownian motion hits the graph of $f_{b_n}$ on the interval $I$, that is a diagonal of the rectangle $I\times f_{b_n}(I)$, and the random variable $Y_{b_n} = \sum_{I \in \mathfrak{C}_{b_n}}\mathbf{1}(Y_n(I))$.

For an interval $I =[x,y] \in \mathfrak{C}_{b_n}$ define $R$ to be the corresponding rectangle $I\times f_b(I)$ and let $R_1$ be the triangle with the vertices $(x,f_b(x))$, $(y,f_b(y))$ and $(x,f_b(x)+2^{-n})$ (so it is the upper left triangle of $R$ with respect to the diagonal of $R$) and $R_2 = R \backslash R_1$ is the lower right triangle part of $R$.

Fix an $n>0$ and let $C$ be the event that $B-f_{b_n}$ has a zero that is contained in $\mathfrak{C}_{b_n}$. Define $\overline{R}$ to be the event that $B-f_b$ has a zero that is contained in $\mathfrak{C}_{b}$ and the corresponding intersection point (by definition it is contained in a rectangle $R$ of described form) of the graph of Brownian motion and the graph of $f_b$ is contained in $R_1$. Analogously, by $\underline{R}$ denote the the event that $B-f_b$ has a zero that is contained in $\mathfrak{C}_{b}$ and the corresponding intersection point of the graph of Brownian motion and the graph of $f_b$ is contained in $R_2$.

Let $\overline{\tau}$ be the first time that $\overline{R}$ happens. Then $\mathbb{P}(B(y)\leq f_{b}(\overline{\tau}))=1/2$. The event $\left\lbrace B(y)\leq f_{b}(\overline{\tau})\right\rbrace $ implies that there is an $s\in I$ such that $B(s)=f_{b_n}(s)$, that is $\mathbb{P}(C |\overline{R}) \geq 1/2$.

Now we go backwards in time. By the time reversal property of Brownian motion the process $\bar{B}(t)= B(2) - B(2-t)$ for $t\in[0,2]$ is again a Brownian motion. Let $\underline{\tau}$ be the first time that $\underline{R}$ happens for the time reversed Brownian motion $B(2-t)$, and let $\bar{x}=2-x$. We want to show that $\mathbb{P}(\bar{B}(\underline{\tau})-\bar{B}(\bar{x})>0|\bar{B}(\underline{\tau})=B(2)-f_{b}(2-\underline{\tau})) \geq \alpha$ for some $\alpha > 0$. In general, for a Brownian motion $B$ the random vector $(B(t)-B(x), B(t))$ has the density
\begin{align*}
\psi(p,q) &= \frac{1}{(2\pi)\sqrt{t(t-x)\frac{x}{t}}}\exp\Big(-\frac{t}{2x}\big(\frac{p^2}{(t-x)} -\frac{2pq\sqrt{t-x}}{t\sqrt{t-x}} + \frac{q^2}{t}\big)\Big) \\
&= \frac{1}{(2\pi)\sqrt{x(t-x)}}\exp\Big(-\frac{1}{2}\big(\frac{tp^2}{x(t-x)} -\frac{2pq}{x} + \frac{q^2}{x}\big)\Big).
\end{align*}
Therefore, with the substitutions $g=B(2)-f_{b}(2-\underline{\tau})$, $p_1 = p - \frac{(\underline{\tau}-\bar{x})g}{\underline{\tau}}$ and $p_2= p_1 \sqrt{\frac{\underline{\tau}}{\bar{x}(\underline{\tau}-\bar{x})}}$ we get
\begin{align*}
 \mathbb{P}(\bar{B}(\underline{\tau})-\bar{B}(\bar{x})>0&|\bar{B}(\underline{\tau})=g)  = \int^{\infty}_0\psi(p,g)dp \\
&\geq\frac{\exp(-\frac{g^2}{2\underline{\tau}})}{(2\pi)\sqrt{\bar{x}(\underline{\tau}-\bar{x})}}\int^{\infty}_{-\frac{\underline{\tau}-\bar{x}}{\underline{\tau}}g}\exp\Big(-\frac{\underline{\tau}p_1^2}{2\bar{x}(\underline{\tau}-\bar{x})}\Big) dp_1 \\
& = \frac{\exp(-\frac{g^2}{2\underline{\tau}})}{(2\pi)\sqrt{\underline{\tau}}}\int^{\infty}_{-\sqrt{\frac{\underline{\tau}-\bar{x}}{\underline{\tau}\bar{x}}}g}\exp\Big(-\frac{p_2^2}{2}\Big)dp_2.
\end{align*}
Since the right hand side is bounded from below, and the event \\$\left\lbrace \bar{B}(\underline{\tau})-\bar{B}(\bar{x})>0 \right\rbrace $ implies the event $C$, we get $\mathbb{P}(C |\underline{R}) \geq \alpha$.
Thus, it follows
\begin{align*}
\mathbb{P}(\zero(B-f_b) \cap C_{b} \neq \emptyset) &\leq \mathbb{P}(\overline{R}\cup\underline{R}) \\
&\leq \mathbb{P}(\overline{R})+\mathbb{P}(\underline{R}) \\
&\leq 2 \mathbb{P}(C |\overline{R})\mathbb{P}(\overline{R}) + \frac{1}{\alpha}\mathbb{P}(C |\underline{R})\mathbb{P}(\underline{R}) \\
&\leq 2 \mathbb{P}(C) + \frac{1}{\alpha}\mathbb{P}(C).
\end{align*}
Therefore,
\begin{align}\label{eq: approximation of probabilities_2}
\mathbb{P}(\zero(B-f_b) \cap C_{b} \neq \emptyset) \leq
(2+\frac{1}{\alpha})
\inf_n
 \mathbb{P}(Y_{b_n}>0).
\end{align}

For a given rectangle that is contained in the square $[1,2]\times [0,1]$ and has the four sides, right side $\textbf{r}$, left side $\textbf{l}$, bottom side $\textbf{b}$ and top side $\textbf{t}$. Call the events that Brownian motion hits these sides $\overline{r},\overline{l},\overline{b}$ and $\overline{t}$, respectively. Let $D$ be the event that Brownian motion hits the diagonal of the rectangle.

Then, analogously to the above argument, by assuming that the events that the graph of Brownian motion hits each side happen instead of the events $\overline{R}$ or $\underline{R}$,
there is a constant $\beta>0$ such that
\begin{equation}\label{eq: approximation of probabilities_3}
\beta \max \left\lbrace \mathbb{P}(\overline{r}), \mathbb{P}(\overline{l}), \mathbb{P}(\overline{b}), \mathbb{P}(\overline{t}) \right\rbrace \leq \mathbb{P}(D).
\end{equation}
If Brownian motion hits the diagonal of the rectangle, then it has to intersect at least one of the sides of the rectangle. That gives the following inequality
\begin{equation}\label{eq: approximation of probabilities_4}
\mathbb{P}(D) \leq 4\max \left\lbrace \mathbb{P}(\overline{r}), \mathbb{P}(\overline{l}), \mathbb{P}(\overline{b}), \mathbb{P}(\overline{t}) \right\rbrace.
\end{equation}

Assume that the graph hits the bottom side of the rectangle and let $\tau'$ be the first such time. Since $\tau$ is a stopping time, by strong Markov property the process $B'(t)=B(t+\tau')-B(\tau')$ is a Brownian motion. If there is a constant $c_1>0$ such that $\sqrt{|\textbf{b}|} \leq c_1|\textbf{r}|$, we can find a constant $\hat{C}>0$, only depending on $c_1$, such that the maximum of $B'$ on the interval $[0,|\textbf{b}|]$ is less than $|\textbf{r}|$ with probability at least $1/\hat{C}$. Since this event implies the event $\overline{r}$, we have $\mathbb{P}(\overline{r}|\overline{b}) \geq 1/\hat{C}$, which implies $\mathbb{P}(\overline{b})\leq \hat{C}\mathbb{P}(\overline{r})$. The inequality $\mathbb{P}(\overline{t})\leq \hat{C}\mathbb{P}(\overline{r})$ can be proven analogously.

If there is a constant $c_2>0$ such that $\sqrt{|\textbf{b}|} \geq c_2 |\textbf{r}|$, then we can show analogously that there is a constant $C'>0$, only depending on $c_2$, such that $\mathbb{P}(\overline{l})\leq C'\mathbb{P}(\overline{t})$, and $\mathbb{P}(\overline{l})\leq C'\mathbb{P}(\overline{b})$.

Note that there are constants $c_3>0$ and $c_4>0$ such that $c_3|\textbf{r}|\leq \mathbb{P}(\overline{r}) \leq c_4|\textbf{r}|$. Assume $\sqrt{|\textbf{b}|} = |\textbf{r}|$. We get that $\frac{c_3}{C'}\sqrt{|\textbf{b}|}\leq \mathbb{P}(\overline{b}) \leq \hat{C}c_4\sqrt{|\textbf{b}|}$.

Therefore, and since the graph of the restriction of $f_{b_n}$ to an interval $I \in \mathcal{C}_{b_n}$ is a diagonal of a rectangle of width $2^{-2n}a_n^2$ and height $2^{-n}$, there are constants $C_1$ and $C_2$ such that
\begin{equation} \label{eq: hitting probabilities for diagonal}
\left.\begin{array}{l l} C_12^{-n},& \text{ if }a_n \leq 1, \\ C_1 2^{-n}a_n,& \text{ if } a_n > 1, \end{array}\right\} \leq \mathbb{P}(Y_n(I)) \leq \left\{\begin{array}{l l} C_22^{-n},& \text{ if }a_n \leq 1, \\ C_22^{-n}a_n,& \text{ if } a_n > 1. \end{array}\right.
\end{equation}

Therefore,
\begin{equation} \label{eq: first moment estimates for diagonal}
\left.\begin{array}{l l} C_1,& \text{ if }a_n \leq 1, \\ C_1a_n,& \text{ if } a_n > 1, \end{array}\right\} \leq
\mathbb{E}(Y_{b_n}) \leq \left\{\begin{array}{l l} C_2,& \text{ if }a_n \leq 1, \\ C_2 a_n,& \text{ if } a_n > 1. \end{array}\right.
\end{equation}

For $I \in \mathfrak{C}_{b_n}$ and $0 \leq \ell < n$, let $I^\ell$ denote the interval from $\mathfrak{C}_{b_\ell}$ that contains $I$, and let $I_1^\ell, I_2^\ell \in \mathfrak{C}_{b_{\ell+1}}$ be the left and right subintervals of $I^\ell$, respectively.

Now we define a binary address $v_1 v_2 \dots v_n$ for $I$. Namely, if $I \subset I_1^\ell$ then $v_{\ell+1}=0$ and if $I \subset I_2^\ell$ then $v_{\ell+1}=1$.

We call an interval $I$ {\em balanced} if the sequence $v_1, \dots, v_n$ contains at least a certain amount of zeros (to be specified later) and otherwise unbalanced. For a balanced interval $I \in \mathfrak{C}_{b_n}$ let $A_I$ denote the event that $I$ is the
leftmost balanced interval for which $Y_n(I)$ happens, and, as before, let $v_1 \dots  v_n$ denote the binary address of $I$. Let $\tau$ be the first time that $f_{b_n}(t)=B(t)$ with $t\in I$ happens.

We take a look again at the event $Z_n(I)$ (see proof of Theorem \ref{thm:cantorset-zeros2}). Fix $n>0$, and an interval $J=[x,y] \in \mathfrak{C}_{b_n}$, so that $I<J$. Let $\widetilde{B}$ be the process
\[\widetilde{B}(t)=(y-\tau)^{-1/2}\Big(B(\tau+(y-\tau) t)-B(\tau)\Big),\]
which is, by the Markov property and Brownian scaling, again a Brownian motion, independent of $\mathcal{F}_{\tau}$. 

Let $\overline{J}$ be the interval $[(y-\tau)^{-1/2}(f_{b_n}(x)-B(\tau)), (y-\tau)^{-1/2}(f_{b_n}(y)-B(\tau))]$ of length $(y-\tau)^{-1/2}2^{-n}$.

Assume that for some $0 \leq \ell < n$ we have $v_{\ell+1}=0$. $I_2^\ell$ contains $2^{n-\ell-1}$ intervals, we label them by $J_1, \dots, J_{2^{n-\ell-1}}$ with $J_j < J_{j+1}$. $\overline{J_j}$ has length at least $(b_\ell)^{-1/2}2^{-n} = \frac{2^{\ell-n}}{a_\ell}$ and at most $(\epsilon b_\ell)^{-1/2}2^{-n}= \frac{2^{\ell-n}}{\sqrt{\epsilon}a_\ell}$ and $\overline{J_j} \subset [(j-1)2^{-n}b_\ell^{-1/2}, (j+2^{n-\ell-1})2^{-n}(\epsilon b_\ell)^{-1/2}]$. Then 

\begin{align}\label{eq: exact conditional probabilities}
\mathbb{P}(Z_n(J_j) \mid A_I) &= \mathbb{P}(\widetilde{B}(1) \in \overline{J_j} \mid A_I)\\
&\geq \int_{(j+2^{n-\ell-1})\frac{2^{\ell-n}}{\sqrt{\epsilon}a_\ell}-\frac{2^{\ell-n}}{a_\ell}}^{(j+2^{n-\ell-1})\frac{2^{\ell-n}}{\sqrt{\epsilon}a_\ell}} \frac{1}{\sqrt{2\pi}}\exp (-\frac{1}{2}t^2) dt.
\end{align}

We can use the following lower bound
\begin{align}\label{eq: conditional probabilities lower estimate}
\mathbb{P}(Z_n(J_j) \mid A_I)  \geq \frac{2^{\ell-n}}{\sqrt{2\pi}a_\ell} \exp\Big(-\frac{1}{2} {((j+2^{n-\ell-1})\frac{2^{\ell-n}}{\sqrt{\epsilon}a_\ell})^2}\Big).
\end{align}
Summing over all $J_1, \dots, J_{2^{n-\ell-1}}$ and over all $\ell$ such that $v_{\ell+1}=0$ gives
\begin{align*}
\mathbb{E} (Z_{b_n} \mid A_I) \geq 
\sum_{1 \leq \ell \leq n: v_\ell=0} \sum_{1\leq j \leq 2^{n-\ell-1}}  \frac{2^{\ell-n}}{\sqrt{2\pi}a_\ell} \exp\Big(-\frac{1}{2} {((j+2^{n-\ell-1})\frac{2^{\ell-n}}{\sqrt{\epsilon}a_\ell})^2}\Big).
\end{align*}
Estimating the inner sum by integration gives
\begin{align}\label{eq: integral estimate}
\sum_{1\leq j \leq 2^{n-\ell-1}}  \frac{2^{\ell-n}}{\sqrt{2\pi}a_\ell} \exp &\Big(-\frac{1}{2} {((j+2^{n-\ell-1})\frac{2^{\ell-n}}{\sqrt{\epsilon}a_\ell})^2}\Big)
\\&\geq \int_1^{2^{n-\ell-1}+1} \frac{2^{\ell-n}}{\sqrt{2\pi}a_\ell} \exp\Big(-\frac{1}{2} {((j+2^{n-\ell-1})\frac{2^{\ell-n}}{\sqrt{\epsilon}a_\ell})^2}\Big) dj\\
&= \sqrt{\frac{\epsilon}{\pi}} \int_{\frac{2^{\ell - n}+1/2}{\sqrt{2\epsilon}a_\ell}}^{\frac{2^{\ell - n}+1}{\sqrt{2\epsilon}a_\ell}} \exp (-t^2) dt.
\end{align}

If $2\epsilon a^2_\ell \leq 2^{\ell - n} + 1/2$ holds we can use the estimate $\int_x^{x+y} \exp(-t^2)dt \geq \frac{\exp{(-x^2)}}{4x}$ (see Figure \ref{fig: integral estimates}).
In case $\epsilon a^2_\ell \leq 1/4$ we will use this estimate and $\int_x^{x+y} \exp(-t^2)dt \geq y \exp (-(x+y)^2)$ otherwise. Using the first estimate we get that (\ref{eq: integral estimate}) is at least $\frac{\sqrt{2}\epsilon a_\ell}{4\sqrt{\pi}(2^{\ell - n}+1/2)}\exp(- (\frac{2^{\ell - n}+1/2}{\sqrt{2\epsilon}a_\ell})^2)$ and in the second $\frac{1}{2\sqrt{2}a_\ell}\exp(- (\frac{2^{\ell - n}+1}{\sqrt{2\epsilon}a_\ell})^2)$.

\begin{figure}

\includegraphics{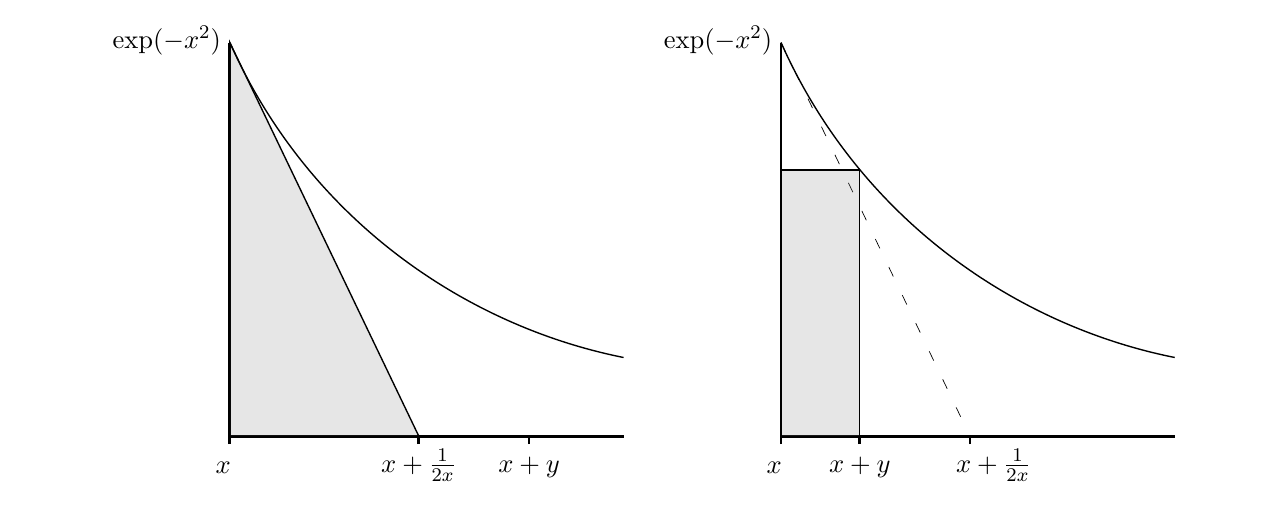}

\caption{If $x+\frac{1}{2x} \leq x+y$, then we bound $\int_x^{x+y} \exp(-t^2)dt$ from below by the area of a triangle. If $x+\frac{1}{2x} > x+y$, then we use the area of a rectangle to bound this integral from below.}

\label{fig: integral estimates}

\end{figure}
Therefore, using the trivial bound $0$ for $\ell=n$ and $\ell=n-1$, we get
\begin{align}\label{eq: balanced first moment}
\mathbb{E} (Z_{b_n} \mid A_I) &\geq \sum_{1 \leq \ell \leq n-2: v_\ell=0, \epsilon a^2_\ell \leq 1/4} \frac{\sqrt{2}\epsilon a_\ell}{4\sqrt{\pi}(2^{\ell - n}+1/2)}\exp(- (\frac{2^{\ell - n}+1/2}{\sqrt{2\epsilon}a_\ell})^2) \\ \nonumber
& \ \ \ \ \  \ \ \ \ \  \ \ \ \ \ + \sum_{1 \leq \ell \leq n-2: v_\ell=0, \epsilon a^2_\ell > 1/4} \frac{1}{2\sqrt{2\pi}a_\ell}\exp(- (\frac{2^{\ell - n}+1}{\sqrt{2\epsilon}a_\ell})^2)
\\&\geq \sum_{1 \leq \ell \leq \frac{5}{6}n: v_\ell=0, \epsilon a^2_\ell \leq 1/4} \frac{\sqrt{2}\epsilon a_\ell}{3\sqrt{\pi}}\exp(- (\frac{2^{\ell - n}+1/2}{\sqrt{2\epsilon}a_\ell})^2) \nonumber
\\& \ \ \ \ \  \ \ \ \ \  \ \ \ \ \ + \sum_{1 \leq \ell \leq n-2: v_\ell=0, \epsilon a^2_\ell > 1/4} \frac{1}{2\sqrt{2\epsilon}a_\ell}\exp(- (\frac{5}{4\sqrt{2\pi}a_\ell})^2). \nonumber
\end{align}

With \eqref{eq: hitting probabilities} we see that
\begin{align}\label{eq: balanced first moment inequal}
1 \geq \Ew (Z_{b_n}) \geq\sum_{I: balanced} \mathbb{E} (Z_{b_n} \mid A_I) \PP(A_I),
\end{align}
where the sum is over all intervals $I$ that are balanced.
That means if we can bound $\mathbb{E} (Z_{b_n} \mid A_I)$ from below by a number $C(n)$ only depending on $n$ (and not depending on the choice of the balanced interval $I$), then $\mathbb{P}(Y_n(I) \text{ for some balanced interval } I)$ is bounded from above by $1/C(n)$.

Now we will consider three cases, namely subsequences $(a_{n_k})$ of $(a_n)$ with $a_{n_k} \to \infty$, $a_{n_k} \to 0$, and otherwise.

We start with the latter case. Here, call 
 an interval $I$ { balanced} if the sequence $v_1, \dots, v_n$ contains at least $n/3$ zeros and otherwise unbalanced. Then, we see that \eqref{eq: balanced first moment} goes to infinity.    

To estimate the probability that $Y_n(I)$ happens for some unbalanced interval $I$ notice that the number of such intervals is bounded from above by $e^{-c_5 n}2^n$ for some constant $c_5>0$. By (\ref{eq: hitting probabilities for diagonal}) this gives
\begin{equation}\label{eq: unbalanced part}
\mathbb{P}(Y_n(I) \text{ for some unbalanced interval } I) \leq   \left\{\begin{array}{l l} e^{-c_5n},& \text{ if }a_n \leq 1, \\ e^{-c_5n} a_n,& \text{ if } a_n > 1. \end{array}\right.
\end{equation}
But note for the case of $a_n > 1$ that we note that $a_n \leq 2^{n/2}$. Thus, (\ref{eq: unbalanced part}) goes to 0 for $n \rightarrow \infty$. It follows that $\inf_n\mathbb{P}(Y_{b_n}>0)=0$.

We proceed with the case that there is a subsequence $(a_{n_k})$ of $(a_n)$ with $a_{n_k} \to 0$. We use the same definition of balanced intervals as in the case before.
Assume there is a sequence $(c_n)$ with $\sum_{n=1}^{\infty} c_n = \infty$ and fixed but arbitrary small $\delta > 0$  and an $n_0>0$ such that for all $n \geq n_0$
\[ \frac{1}{c_n}\geq a_n \geq \sqrt{\frac{1}{(8-\delta)\epsilon\ln{\frac{1}{c_n}}}}. \]
Now take the subsequence $(a_{n_k})$ such that for all $n_k \geq n_0$
\[ \frac{1}{2\sqrt{\epsilon}}\geq a_{n_k} \geq \sqrt{\frac{1}{(8-\delta)\epsilon\ln{\frac{1}{c_{n_k}}}}}. \]
The right inequality is equivalent to
\[ -\frac{1}{(8-\delta)\epsilon a_{n_k}^2}  \geq \ln{c_{n_k}}. \]
From this it follows that
\[ -\frac{1}{8\epsilon a_{n_k}^2} + \ln{a_{n_k}}  \geq \ln{c_{n_k}}, \]
and
\[ \exp{(-\frac{1}{8\epsilon a_{n_k}^2})} \cdot {a_{n_k}}  \geq {c_{n_k}}. \]
With (\ref{eq: balanced first moment}) and the argument following (\ref{eq: balanced first moment}) we have $$\inf_n\mathbb{P}(Y_n(I) \text{ for some unbalanced interval } I)=0.$$
Thus together with (\ref{eq: unbalanced part}) it follows $\inf_n\mathbb{P}(Y_{b_n}>0)=0$.

To finish we look at the case of having a subsequences $(a_{n_k})$ of $(a_n)$ with $a_{n_k} \to \infty$. Let $d_n := \frac 1 2 \sum_{i=1}^{n} \frac 1 {a_i}$. Now we define an interval $I$ to balanced if the corresponding binary address $v_1, \dots, v_n \in\{0,1\}^n$ fulfills that $\sum_{i=1}^{n} \frac {v_i} {a_i} \leq d_n$.

Observe that we can apply the exceptional Chebychev inequality, for a positive number $k$,
\begin{align}\label{chebychev ineq}
\PP \big( \sum_{i=1}^{n} \frac {v_i} {a_i} \geq d_n \big) &\leq \Ew \Big( \exp \big[k(\sum_{i=1}^{n} \frac {v_i} {a_i} - d_n) \big] \Big)\\ \nonumber
&= \exp(-kd_n)\cdot \prod_{i=1}^{n} \Ew \big( \exp (k \frac {v_i} {a_i}) \big)\\ \nonumber
&= \exp(-kd_n) \frac 1 2 \big(1 + \exp ( \frac k {a_i}) \big).
\end{align}

Note that if $\sum_{i=1}^{n} \frac 1 {a_i} \leq C \log n$, then we can find a $k>0$ for \eqref{chebychev ineq} such that $\PP \big( \sum_{i=1}^{n} \frac {v_i} {a_i} \geq d_n \big) \cdot a_n \to 0$. By (\ref{eq: hitting probabilities for diagonal}) it follows that the probability that $Y_n(I)$ happens for some unbalanced interval
goes to $0$,
 and also $\mathbb{P}(Y_n(I) \text{ for some balanced interval } I)$ goes to $0$ as $n \to \infty$ by (\ref{eq: balanced first moment}) and the argument following (\ref{eq: balanced first moment}).

The claim follows now from (\ref{eq: approximation of probabilities_2}).

\end{proof}

If we require the sequence $(a_n)$ to fulfill $a_n^2 - \frac{1}{2}a_{n+1}^2 \geq x_n$ for some positive sequence $(x_n)$ with $x_n< a_n^2$ for all $n$ instead of the condition $a_n^2 - \frac{1}{2}a_{n+1}^2 \geq \epsilon a_n^2$ for all $n$ that we used so far, then the analogue to the Theorems \ref{thm:cantorset-zeros2} and \ref{thm: cantorset-zeros 2nd part} is the following result.

\begin{corollary}\label{rm: condition on a_n's}
If the sequence $(a_n)$ fulfills $a_n^2 - \frac{1}{2}a_{n+1}^2 \geq x_n$ for some positive sequence $(x_n)$ with $x_n< a_n^2$ for all $n$, then
$\mathbb{P}(\zero(B-f_b) \cap C_b \neq \emptyset)>0$ holds if either $\sum_{n=1}^{\infty} \frac{a^2_n}{\sqrt{x_n}} < \infty$, or $\sum_{n=1}^{\infty} \frac{1}{\sqrt{x_n}} < \infty$, and

$\mathbb{P}(\zero(B-f_b) \cap C_b \neq \emptyset)=0$ holds if

(i) there is a sequence $(c_n)$ with $\sum_{n=1}^{\infty} c_n = \infty$ and a fixed but arbitrary small $\delta > 0$ and an $n_0>0$ such that for all $n \geq n_0$
\[ \frac{1}{c^2_n} \geq x_n \geq \frac{1}{(8-\delta)\ln{\frac{1}{c_n}}}, \]
or if

(ii)$\sum_{\ell = 1}^n \frac{1}{\sqrt{x_\ell}} \rightarrow \infty$ for $n \rightarrow \infty$ and if there is a constant $C>0$ such that $\sum_{\ell = 1}^n \frac{1}{\sqrt{x_\ell}} \leq C \log n$ for all $n>0$.
\end{corollary}

\begin{proof}
Analogously to the proofs of the Theorems \ref{thm:cantorset-zeros2} and \ref{thm: cantorset-zeros 2nd part}.
\end{proof}

Theorems \ref{thm:cantorset-zeros2} and \ref{thm: cantorset-zeros 2nd part} do not give an answer for certain sequences $(a_n)$, for instance $a_n=\frac{1}{n}$, whether or not the zero set of $B-f_b$ contains points of the corresponding generalized Cantor set with positive probability. It is natural to ask if we can strengthen the methods we used to get a stronger result.

Note that by (\ref{eq: approximation of probabilities_3}) and (\ref{eq: approximation of probabilities_4}) the probability of the event that Brownian motion hits a diagonal of a rectangle is up to a constant the maximum of the probabilities of the events that that Brownian motion hits the lower horizontal side or the right vertical side of the rectangle. Thus, it worth looking at the event that Brownian motion hits the lower horizontal side of the rectangle instead of the event $Z_n(J_j)$ for the estimate (\ref{eq: exact conditional probabilities}). Call the former event $X_n(J_j)$, then we get instead of (\ref{eq: exact conditional probabilities})

\begin{align*}
\mathbb{P}(X_n(J_j) \mid A_I)  \geq \int_{x}^{x+ 2^{-2n} a_n^2} \frac{1}{\sqrt{2\pi t}}\exp (-\frac{{(2^{-\ell-1} + 2^{-n}j)}^2}{2t}) dt
\end{align*}

for some point $x \in [\sqrt{\epsilon}2^{-2\ell} a_\ell^2, 2^{-2\ell} a_\ell^2)$. But here we see that summing over all $j$ and $\ell$ would not give an expression that goes to infinity for $n \rightarrow \infty$ for any possible sequence $(a_n)$. Thus, looking at the event that that Brownian motion hits the lower horizontal side cannot provide an improvement of the result \ref{thm: cantorset-zeros 2nd part}.

Comparing
\begin{multline*}
\int_{(j+2^{n-\ell-1})\frac{2^{\ell-n}}{\sqrt{\epsilon}a_\ell}-\frac{2^{\ell-n}}{a_\ell}}^{(j+2^{n-\ell-1})\frac{2^{\ell-n}}{\sqrt{\epsilon}a_\ell}} \frac{1}{\sqrt{2\pi}}\exp (-\frac{1}{2}t^2) dt \\ \leq
\frac{2^{\ell-n}}{\sqrt{2\pi}a_\ell} \exp\Big(-\frac{1}{2} {((j+2^{n-\ell-1})\frac{2^{\ell-n}}{\sqrt{\epsilon}a_\ell}-\frac{2^{\ell-n}}{a_\ell})^2}\Big),
\end{multline*}
to the estimate (\ref{eq: conditional probabilities lower estimate}), and
\begin{align*}
\sum_{1\leq j \leq 2^{n-\ell-1}}& \frac{2^{\ell-n}}{\sqrt{2\pi}a_\ell} \exp\Big(-\frac{1}{2} {((j+2^{n-\ell-1})\frac{2^{\ell-n}}{\sqrt{\epsilon}a_\ell}-\frac{2^{\ell-n}}{a_\ell})^2}\Big) \\ &\leq
\int_0^{2^{n-\ell-1}} \frac{2^{\ell-n}}{\sqrt{2\pi}a_\ell} \exp\Big(-\frac{1}{2} {((j+2^{n-\ell-1})\frac{2^{\ell-n}}{\sqrt{\epsilon}a_\ell}-\frac{2^{\ell-n}}{a_\ell})^2}\Big) dj \\
&= \sqrt{\frac{\epsilon}{\pi}} \int_{\frac{\frac{1}{2}-\sqrt{\epsilon}2^{\ell - n}}{\sqrt{2\epsilon}a_\ell}}^{\frac{1-\sqrt{\epsilon}2^{\ell - n}}{\sqrt{2\epsilon}a_\ell}} \exp (-t^2) dt \\
&\leq \frac{1}{2\sqrt{2\epsilon}a_\ell}\exp(- (\frac{\frac{1}{2}-\sqrt{\epsilon}2^{\ell - n}}{\sqrt{2\epsilon}a_\ell})^2),
\end{align*}
to (\ref{eq: integral estimate}) and (\ref{eq: balanced first moment})
we see that
for the inequality (\ref{inequality of thm: cantorset-zeros 2nd part}) we cannot improve the result \ref{thm: cantorset-zeros 2nd part} by using better bounds.

For the proof of Theorem \ref{thm:cantorset-zeros2} it is easy to check that better estimates can not provide a stronger result. In particular, by bounding $S(z)$ from above we see that the second moment can at most be improved by a constant factor.


\begin{remark}
Theorems \ref{thm:cantorset-zeros2} and \ref{thm: cantorset-zeros 2nd part} can be extended to more classes of Cantor-like functions. We consider for a given positive, real sequence $(a_n)$ a sequence $(b_{k,n})$ defined by $b_{k,n}:= (k^{-n} \cdot a_{n})^2$ with $k>0$. For this sequence we define a corresponding $k$-Cantor-type set and denote it by $C_{b_k}$.
Take a closed interval $I$ of length $b_{k,0}$. Define $\mathfrak{C}_{b_{k,1}}$ as the set consisting of $k$ disjoint closed subintervals of $I$ of length $b_{k,1}$, the left one (for which the left endpoint coincides with the left endpoint of $I$), the right one (for which the right endpoint coincides with the right endpoint of $I$) and $k-2$ intervals such that two neighboring intervals have the distance $\frac{b_{k,0} - kb_{k,1}}{k-1}$.
Continue recursively, if $J \in \mathfrak{C}_{b_{n}}$, then include in the set $\mathfrak{C}_{b_{n +1}}$ its $k$ closed subintervals of length $b_{k,n+1}$. Define the set $C_{b_{n}}$ as the union of all the intervals from $\mathfrak{C}_{b_{n}}$.
For any $n$, the family $\mathfrak{C}_{b_{k,n}}$ is the set of all connected components of the set $C_{b_{k,n}}$.
The $k$-Cantor set is a compact set defined as $C_{b_k}=\bigcap_{n=1}^\infty C_{b_{k,n}}$.
Now we construct a $k$-Cantor-type function corresponding to the $k$-Cantor type set above. Define the function $f_{b_{k,1}}$ so that it has values $0$ and $1$ at the left and the right endpoint of the interval $I$, respectively, values $1/k$ on the most left, $2/k$ on the second most left, ..., and $(k-1)/k$ on the least most left of the $k-1$ disjoint intervals of the set $I\backslash {C}_{b_{k,1}}$, and interpolate linearly on the intervals in $\mathfrak{C}_{b_{k,1}}$. Recursively, construct the function $f_{b_{k,n +1}}$ so that for every interval $J =[s,t]\in \mathfrak{C}_{b_{k,n}}$, the function $f_{b_{k,n +1}}$ agrees with $f_{b_{k,n}}$ at $s$ and $t$, it has values $(f_{b_{k,n}}(t)-f_{b_k,n}(s))\cdot\frac{1}{k} +f_{b_k,n}(s)$ on the most left, $(f_{b_{k,n}}(t)-f_{b_k,n}(s))\cdot\frac{2}{k} +f_{b_k,n}(s)$ on the second most left, ..., and $(f_{b_{k,n}}(t)-f_{b_k,n}(s))\cdot\frac{k-1}{k} +f_{b_k,n}(s)$ on the least most left of the $k-1$ disjoint intervals of the set $J \backslash C_{b_{k,n +1}}$ and interpolate linearly on the intervals in $\mathfrak{C}_{b_{k,n +1}}$. See Figure \ref{fig:triplecantor} for an example of $k=3$.

\begin{figure}
\includegraphics{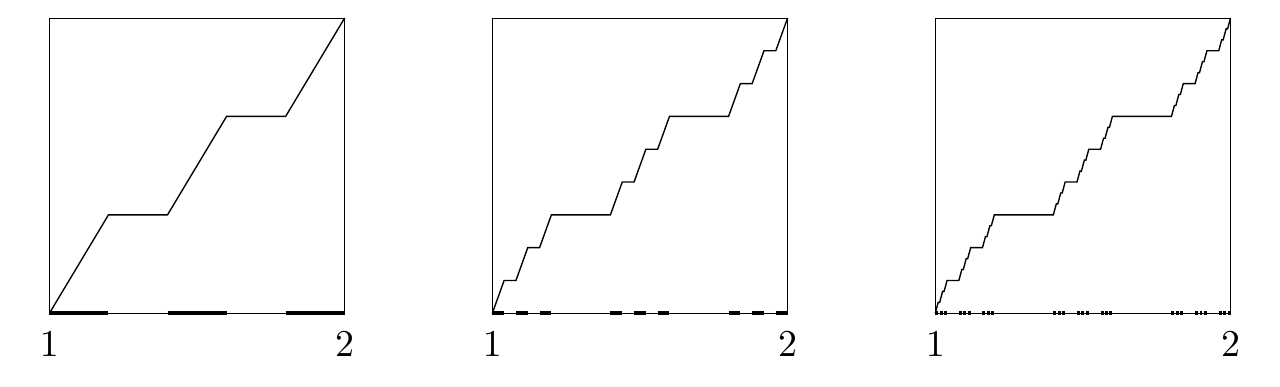}

\caption{First three approximations of the $3$-Cantor type function on the interval $[1,2]$ (functions $f_{b_{3,n}}$ from the construction for $b_1=0.2$, $b_2=0.04$ and $b_3=0.01$, or $a_1=\sqrt{1.8}$, $a_2=1.8$ and $a_3=2.7$, respectively. Approximations of the Cantor set $C_{b_3}$ are drawn in bold).}

\label{fig:triplecantor}

\end{figure}

As in section \ref{section: Cantor function def}, for a fixed $k>1$ we fix an arbitrary small $\epsilon > 0$ and require the sequence $(a_{n})$ to fulfill the condition $a_{n}^2 - \frac{1}{k}a_{n+1}^2 \geq \epsilon a_{n}^2$ (or equivalently $b_{k,n}- kb_{k,n+1} \geq \epsilon b_{k,n}$) for all $n$. Then, Theorems \ref{thm:cantorset-zeros2} and \ref{thm: cantorset-zeros 2nd part} hold also for the $k$-Cantor-type function.

Further, note that if we have the weaker condition $b_{k,n}- kb_{k,n+1} > 0$ for all $n$, then we can apply Corollary \ref{rm: condition on a_n's}.

\end{remark}

\section{Isolated zeros - general criteria}\label{Isolated zeros - general results}

We will now state two criteria determining whether a zero of $B-f$ for a continuous function $f$ is almost surely isolated or not isolated.
For any function $g$ defined on some subset (or the whole) of $\mathbb{R}^+$ denote by $\zero(g)$ the set of zeros of $g$ in $(0,\infty)$. Recall that we denoted the middle $\alpha$-Cantor function by $f_\alpha$.
Antunovi\'{c}, Burdzy, Peres, and Ruscher showed that for every $\alpha<1/2$ there is an $\alpha$-\holder continuous function $f \colon \mathbb{R}^+ \to \mathbb{R}$ such that the set $\zero(B-f)$ has isolated points with positive probability, see Theorem 1.2 of \cite{ABPR}.
To show this result they proved that the zero set of $B-f_\alpha$ has isolated points with positive probability for $\alpha > 1/2$. We want to extend that result for our generalized
class of Cantor functions in the next section. First we need to look at a criterion for having isolated zeros.

\begin{proposition}\label{prop:basic_facts2}

Let $f \colon \mathbb{R}^+ \to \mathbb{R}$ be a continuous function.

\begin{itemize}

\item[(i)]

Let $A$ be a closed subset of $\mathbb{R}^+$ such that for any $t \in A$
\[\liminf_{s \to t}\frac{|f(s)-f(t)|}{\sqrt{2|t-s|\log\log{\frac{1}{|t-s|}}}} > 1. \]
Then, almost surely any point in $\zero(B-f) \cap A$ is isolated in $\zero(B-f)$.

\item[(ii)]

Let $A \subset \mathbb{R}^+$ be a set such that for any $t \in A$
\[\limsup_{s \to t}\frac{|f(s)-f(t)|}{\sqrt{2|t-s|\log\log{\frac{1}{|t-s|}}}} < 1. \]
Then, almost surely any point in $\zero(B-f) \cap A$ is not isolated in $\zero(B-f)$.

\end{itemize}

\end{proposition}

\begin{proof}
(i)
We define a sequence of stopping times $(\tau_n)$. Let $$\tau_0 = \min\left\{ t\in A : B(t)=f(t)\right\},$$ and $$\tau_n = \min\left\{ t\in A, t> \tau_{n-1}: B(t)=f(t)\right\}.$$
Since $\tau_n$ is a stopping time for every $n$, $t \mapsto B(\tau_n+t)-B(\tau_n)$ is a Brownian motion if $\tau_n < \infty$. We can apply the law of iterated logarithm (see Theorem 5.1 in \cite{MP}) to get that almost surely for all $n$ we have

\[
\liminf_{t \downarrow 0}\frac{B(\tau_n+t)-B(\tau_n)}{\sqrt{2t\log\log{\frac{1}{t}}}} = - 1 \ \text{ and } \limsup_{t \downarrow 0}\frac{B(\tau_n+t)-B(\tau_n)}{\sqrt{2t\log\log{\frac{1}{t}}}} = 1.
\]

Thus all $\tau_n$'s are isolated from the right. By the reverse property of Brownian motion all $\tau_n$'s are also isolated from the left. This implies that $\tau_n$ converges to $\infty$ since A is a closed set. Therefore, every zero in $A$ is contained in the sequence $(\tau_n)$.

(ii) Assume that there exists a an isolated zero in the set $A$ with positive probability. Then, there is a $q \in \mathbb{Q}$ such that $\tau_q = \min\left\{ t\geq q : B(t)=f(t)\right\}$ is an isolated zero in $A$. 

Since $\tau_q$ is a stopping time, the process $B_q(t)=B(\tau_q+t)-B(\tau_q)$, is, by the strong Markov property, a Brownian motion independent of the sigma algebra $\mathcal{F}_{\tau_q}$. By the law of the iterated logarithm it follows that $\tau_q$ is not isolated.

\end{proof}

Note that (i) is a stronger statement than Proposition 2.2.(i) of \cite{ABPR} for closed subsets of $\mathbb{R}^+$. According to Proposition 2.2.(ii) of \cite{ABPR} almost surely all isolated points of $\zero(B-f)$ are located inside the set $A_f^+ \cup A_f^-$, where

$A_f^+ = \{t\in \mathbb{R}^+: \lim_{h \downarrow 0}\frac{f(t+h)-f(t)}{\sqrt{h}} = \infty\}$ and

$A_f^- = \{t\in \mathbb{R}^+: \lim_{h \downarrow 0}\frac{f(t+h)-f(t)}{\sqrt{h}} = -\infty\}$. In particular, (ii) shows that not all points in $\lbrace A_f^+ \cup A_f^-\rbrace \cap \zero(B-f) $ have to be isolated in $\zero(B-f)$.


\section{Isolated zeros of Brownian Motion minus Cantor function}

In \cite{ABPR} it was shown for the middle $\alpha$-Cantor function $f_\alpha$ that $B-f_\alpha$ has isolated zeros with positive probability if $\alpha > 1/2$. Recall that for $a_n = \frac{1}{x^{n}}$ with some $x>1$ gives $b_{n}:= (2^{-n} \cdot \frac{1}{x^{n}})^2 = (2x)^{-2n}$ which corresponds to $\alpha = 1- \frac{1}{x}$. Therefore, the following proposition extends this result of \cite{ABPR}.

\begin{proposition}
The set $\zero(B-f_b)$ has isolated points with positive probability if $\sum_{n=1}^{\infty} a_n < \infty$, and no isolated points almost surely if $\sum_{n=1}^{\infty} \frac{1}{a_n} < \infty$.
\end{proposition}

\begin{proof}
First we will prove that if $\sum_{n=1}^{\infty} a_n < \infty$, then $\zero(B-f_b)$ has isolated points with positive probability.
For $A \subset \mathbb{R}^+$ define $Z(A)$ as the event $\{\zero(B-f_b) \cap A \neq \emptyset\}$.
We claim that there exists a constant $c_1$, such that for any interval $J \subset [0,1]$ of length $|J|$, we have

\begin{equation}\label{eq:cuts_probabilities2}
\mathbb{P}(Z(C_b \cap f_b^{-1}(J))) \leq c_1 |J|.
\end{equation}

In order to show this statement, fix an interval $J$ and take the biggest integer $n$ satisfying $|J| \leq 2^{-n}$. Notice that $J$ can be covered by two consecutive binary intervals $J_1$ and $J_2$ of length $2^{-n}$.
Moreover, there are consecutive $I_1,I_2 \in \mathfrak{C}_{b_n}$ such that $f_b(I_i) = J_i$ for  $i=1,2$,
and $C_b \cap  f_b^{-1}(J) \subset I_1 \cup I_2$. 

Now we will again use the notation that we introduced in the proof of Theorem \ref{thm:cantorset-zeros2}.
Assume that $\zero(B-f_b)\cap C_b \neq \emptyset$ and denote the first zero of $B(t)-f_b(t)$ in the generalized Cantor set $C_b$ by $\tau$ ($\tau$ exists since $\zero(B-f_b)\cap C_b$ is a closed set). For an interval $I = [s,t] \in \mathfrak{C}_{b_n}$ assume that $\tau \in I$. Since $\tau$ is a stopping time, and by Brownian scaling, the conditional probability $\mathbb{P}\Big(Z_n(I) \mid \mathcal{F}_\tau, \tau \in I\Big)$ is equal to the probability that Brownian motion at time $1$ is between $y_1=(f_b(s)-f_b(\tau))(t-\tau)^{-1/2}$ and $y_2=(f_b(t)-f_b(\tau))(t-\tau)^{-1/2}$. Since $f_b(s) \leq f_b(\tau) \leq f_b(t)$ we see that $y_1 \leq 0$ and $y_2 \geq 0$. Moreover, $t-\tau \leq 4^{-n} = (f_b(t)-f_b(s))^2$ leads to $y_2-y_1 \geq 1$.

Thus we can bound the probability
\begin{equation}\label{eq: approximation of probabilities_1-5}
\mathbb{P}\Big(Z_n(I)
\mid
 \mathcal{F}_\tau, \tau \in I \Big) 
\geq \mathbb{P}(0 \leq  B(1) \leq 1) = \alpha,
\end{equation}
 for some $\alpha>0$. 
Hence,
\begin{equation} \label{eq:approximation_by_diagonals_for_small_gamma}
{\alpha}\mathbb{P}(Z(C_b \cap I_i)) \leq  \mathbb{P}(Z_n(I_i)).
\end{equation}

But by the first inequality in (\ref{eq: hitting probabilities}), the probability on the right hand side of (\ref{eq:approximation_by_diagonals_for_small_gamma}) is bounded from above by $2^{-n}$. Applying this fact in (\ref{eq:approximation_by_diagonals_for_small_gamma}) and summing the expression for $i=1,2$, we obtain (\ref{eq:cuts_probabilities2}).

By Theorem \ref{thm:cantorset-zeros2} the set $\zero(B-f_b) \cap C_b$ is non-empty with some probability $p>0$.
For every $n>0$ take an arbitrary $b'_n=(2^{-n}a'_n)^2$ such that $b_n < b'_n < 2^{-2n}$ and $n_0$ such that $\sum_{n\geq n_0} a'_n \leq p/(2c_1)$.

For $n \geq n_0$ consider the interval $J_{k,n}=[k2^{-n}-\frac{1}{2}\sqrt{b'_n}, k2^{-n}+\frac{1}{2}\sqrt{b'_n}]$ and define the set $M_{n_0}= \bigcup_{n \geq n_0} \bigcup_{0 \leq k \leq 2^{n}}J_{k,n}$. By (\ref{eq:cuts_probabilities2}) and the choice of $n_0$, we see that $\mathbb{P}(Z(C_b \cap f_b^{-1}(M_{n_0}))) \leq p/2$. Hence, the event that there is a zero of $B(t)-f_{b}(t)$ in the set $C_b \cap \text{Int}(C_{b_{n_0}})\backslash f_b^{-1}(M_{n_0})$ has probability of at least $p/2$ (here $\text{Int}(C_{b_{n_0}})$ is the interior of the set $C_{b_{n_0}}$). Now the claim follows if we prove that any such zero is isolated. Take $t \in C_b \cap \text{Int}(C_{b_{n_0}})\backslash f_b^{-1}(M_{n_0})$ and any $s \ne t$ in the same connected component of $\text{Int}(C_{b_{n_0}})$. The biggest integer $\ell$ such that both $s$ and $t$ are contained in the same interval of $\mathfrak{C}_{b_\ell}$ satisfies $\ell \geq n_0$.
Further, $|f_b(s)-f_b(t)| \geq \frac{1}{2}\sqrt{b'_\ell}$ and $b_{\ell+1} \leq |s-t| \leq b_{\ell}$. With Proposition \ref{prop:basic_facts2}(i) the claim follows by taking for example $a'_n = \sqrt{a_n}$ for $n$ big enough.

For the second part of the claim we just need to apply the Proposition 2.2.(i) of \cite{ABPR} (see above).
\end{proof}

\section*{Acknowledgments}
The author thanks gratefully Michael Scheutzow for fruitful discussions and advice.


\begin{thebibliography}{4}



\bibitem[ABPR]{ABPR}
T.~Antunovi\'{c}, K.~Burdzy, Y.~Peres, and J.~Ruscher.
\newblock {\em Isolated zeros for Brownian motion with variable drift}.
\newblock Electronic Journal of Probability. Vol. 16, No. 65: 1793--1814, 2011.


\bibitem[APV]{APV}
T.~Antunovi\'{c}, Y.~Peres, and B.~Vermesi.
\newblock {\em Brownian motion with variable drift can be space-filling}.
\newblock Proc. Amer. Math. Soc. Vol. 139: 3359-3373, 2011.


\bibitem[HT]{HT}
G. G. Hamedani, M. N. Tata,
\newblock {\em On the determination of the bivariate normal distribution from distributions of linear combinations of the variables}.
\newblock The American Mathematical Monthly. Vol. 82, No. 9: 913--915, 1975.


\bibitem[MP]{MP}
P.~M\"{o}rters and Y.~Peres.
\newblock {\em Brownian Motion}.
\newblock Cambridge Series in Statistical and Probabilistic Mathematics.
  Cambridge University Press, 2010.


\bibitem[PS]{PS}
Y. Peres and P. Sousi.
\newblock {\em Brownian motion with variable drift: 0-1 laws, hitting probabilities and Hausdorff dimension.}
\newblock Available at http://arxiv.org/abs/1010.2987v1.


\bibitem[R12]{R12}
J. Ruscher.
\newblock {\em A note on fast times of Brownian motion with variable drift.}
\newblock Preprint, 2012.


\bibitem[RY]{RY}
Daniel Revuz and Marc Yor.
\newblock {\em Continuous martingales and {B}rownian motion}, volume 293 of
  {\em Grundlehren der Mathematischen Wissenschaften [Fundamental Principles of
  Mathematical Sciences]}.
\newblock Springer-Verlag, Berlin, third edition, 1999.


\bibitem[TW]{TaylorWatson}
S.~J. Taylor and N.~A. Watson.
\newblock {\em A {H}ausdorff measure classification of polar sets for the heat
  equation}.
\newblock {Math. Proc. Cambridge Philos. Soc.}, Vol. 97, No. 2: 325--344, 1985.


\end{thebibliography}
\end{document}